\documentclass[12pt]{article}
\usepackage{geometry,latexsym,amssymb,verbatim,cite}
\newtheorem{theorem}{Theorem}

\newtheorem{lemma}[theorem]{Lemma}
\newtheorem{corollary}[theorem]{Corollary}

\begin{document}


\title{\bf Finite Sholander Trees, Trees,\\ and their Betweennes}
\author{
{\bf Va\v{s}ek Chv\'{a}tal}$^1$,
{\bf Dieter Rautenbach}$^2$, 
and {\bf Philipp Matthias Sch\"{a}fer}$^2$}
\date{}
\maketitle

\begin{center}
$^1$
Canada Research Chair in Combinatorial Optimization, Department of Computer Science and Software Engineering, 
Concordia University, 
Montr\'{e}al, 
Qu\'{e}bec,
Canada,
email: \texttt{chvatal@cse.concordia.ca}\\[3mm]
$^2$
Institut f\"{u}r Optimierung und Operations Research, 
Universit\"{a}t Ulm, 
D-89069 Ulm,
Germany,
emails: \texttt{dieter.rautenbach@uni-ulm.de, philipp.schaefer@gmail.com}
\end{center}

\begin{abstract}
  We provide a proof of Sholander's claim (Trees, lattices, order, and
  betweenness, {\it Proc. Amer. Math. Soc.} {\bf 3}, 369-381 (1952))
  concerning the representability of collections of so-called segments
  by trees, which yields a characterization of the interval function of
  a tree.  Furthermore, we streamline Burigana's characterization (Tree
  representations of betweenness relations defined by intersection and
  inclusion, {\it Mathematics and Social Sciences} {\bf 185}, 5-36
  (2009)) of tree betweenness and provide a relatively short proof.
\end{abstract}

{\small {\bf Keywords:} Graph; tree; interval function; betweenness; convexity

{\small {\bf MSC 2010 classification: 05C05, 52A01, 52A37}}

\section{Introduction}

Trees form one of the most simple yet important classes of graphs
with countless applications ranging from data structures and VLSI
design over mathematical psychology to gardening.  Here we consider
two closely related papers on trees: one by Sholander~\cite{sh},
published in 1952, and the other by Burigana~\cite{bu}, published in
2009.  Both of these papers include characterizations of certain
ternary relations associated with trees.

We use the term {\it tree\/} in the sense defined by K\H{o}nig
($\!\!$\cite{ko}, p.47): a finite, simple, undirected, and connected
graph without cycles. Sholander~\cite{sh} used this term in a
different sense: he studied collections of so-called {\it segments},
which are subsets of a set $V$ indexed by all ordered pairs of
elements of $V$ and he referred to such a collection as a {\it tree}
if it satisfies certain postulates. He stated without a proof that
these postulates characterize the function that assigns to every pair
of vertices of a tree in the sense of K\H{o}nig the set of vertices on
the path joining these two vertices; nowadays, this function is called
the {\it interval function} of the tree~\cite{mune}. Interval
functions of K\H{o}nig trees are easily seen to be trees in
Sholander's sense, but it is not obvious that all finite Sholander trees are
representable as interval functions of K\H{o}nig trees. In Section 2,
we supply the missing proof of this claim. 

The {\it tree betweenness\/} of a tree $T$} is defined as the set of
all ordered triples $(x,y,z)$ such that $x,y,z$ are (not necessarily
distinct) vertices of $T$ and $y$ belongs to the path in $T$ that
joins $x$ and $z$; the {\it strict tree betweenness\/} of $T$ is
defined as the set of all ordered triples $(x,y,z)$ such that $x,y,z$
are pairwise distinct vertices of $T$ and $y$ belongs to the path in
$T$ that joins $x$ and $z$. It is a routine matter to restate
Sholander's characterization of the interval function of a tree as a
characterization of tree betweenness; this was done, with refinements,
by Sholander himself in the same paper~\cite{sh}; subsequently, Defays
\cite{de} found another characterization of tree betweenness. Burigana
(Theorem 1 in~\cite{bu}) characterized strict tree betweenness by a
list of five properties that do not involve the notion of a tree. His
proof is spread over some seven pages; in Section 3, we give a shorter
proof; actually, we prove a simpler theorem, of which Burigana's is an
instant corollary.  In addition, we restate the simplified
characterization of strict tree betweenness in terms of tree
betweenness.

Before proceeding to our results, let us put their subject in a
broader context by mentioning a few related references.  Mulder and
Nebesk\'{y}~\cite{mune,ne1,ne2} studied interval functions of
arbitrary graphs. Tree betwenness is a special kind of {\em metric
  betweenness\/} that, for a prescribed metric space $(V,{\rm dist})$,
consists of all ordered triples $(x,y,z)$ such that $x,y,z$ are (not
necessarily distinct) points of $V$ and ${\rm dist}(x,y)+{\rm
  dist}(y,z)={\rm dist}(x,z)$. This concept was first studied by
Menger~\cite{me} in 1928; references to subsequent work on it can be
found in~\cite{ch1}. Another special kind of metric betweeness is {\em
  Euclidean betweenness,\/} where the metric space is a Euclidean
space or some subspace of it. In his development of geometry, Euclid
used the notion of betweenness only implicitly; its explicit
axiomatization was first carried out by Pasch \cite{Pas} and then
gradually refined by Peano \cite{Pea}, Hilbert \cite{Hil}, Veblen
\cite{Veb}, and Huntington and Kline~\cite{hukl}. In particular,
Huntington and Kline suggested the study of other ternary relations
(meaning subsets of $V^3$, where $V$ is some set) that resemble
Euclidean betweenness: for example, they mention the set of all
ordered triples $(x,y,z)$ such that $x,y,z$ are natural numbers and
$y=xz$. Pitcher and Smiley~\cite{PS} continued in this direction.
Another particular kind of betwenness is {\em order betweenness\/}
that, for a prescribed partially ordered set $(V,\preceq)$, consists
of all ordered triples $(x,y,z)$ such that $x,y,z$ are (not
necessarily distinct) points of $V$ and $x\preceq y\preceq z$ or
$z\preceq y\preceq x$.  This concept was first studied by
Birkhoff~\cite{bi} in 1948; Altwegg~\cite{al} characterized order 
betweenness by a list of six properties that do not involve the notion
of a partially ordered set; subsequently, Sholander~\cite{sh} and
D\"{u}ntsch and Urquhart \cite{duur} found other characterizations of
order betweenness.

\section{Finite Sholander Trees are Trees}

Sholander studies mappings that assign to each ordered pair $(a,b)$ of
elements of a set $V$ a subset of $V$, which we denote as $[ab]$. From
postulates 
\begin{center}
\begin{tabular}{lp{14.5cm}}
(S) & $\forall a,b,c\in V : \exists  d\in V: [ab]\cap [bc]=[bd]$,\\
(T) & $\forall a,b,c\in V : [ab]\subseteq [ac] \:\Rightarrow\: [ab]\cap [bc]=\{ b\}$,
\end{tabular}
\end{center}
he derives a number of corollaries that include
\begin{center}
\begin{tabular}{lp{14.5cm}}
(1.2) & $\forall a,b\in V : b\in [ab]$,\\
(1.4) & $\forall a,b\in V : [ab]=[ba]$,\\
(1.5) & $\forall a,b,c\in V : b\in [ac]\Leftrightarrow [ab]\subseteq [ac]$,\\
(1.7) & $\forall a,b,c\in V : (b\in [ac]\wedge c\in [ab])\:\Rightarrow\: b=c$,\\
(1.10) & $\forall a,b,c,d\in V : [ab]\cap [bc]=[bd]\:\Rightarrow\: [ad]\cap [dc]=\{ d\}$.
\end{tabular}
\end{center}
(The labels (S), (T), (1.2), etc. used in this section are copied
directly from \cite{sh} for ease of reference.) Then he defines a
tree as a mapping $(u,v)\mapsto [uv]$ from $V^2$ to $2^V$ that
satisfies (S), (T), and
\begin{center}
\begin{tabular}{lp{14.5cm}}
(U$_1$) & $\forall a,b,c\in V : [ab]\cap [bc]=\{ b\}\:\Rightarrow\: [ab]\cup [bc]=[ac]$.
\end{tabular}
\end{center}
Having noted that K\H{o}nig defined a tree as a finite connected graph
that contains no cycles, he states (\cite{sh}, p. 370) that ``{\it
  Trees in our sense which are finite are trees in K\H{o}nig's
  sense.}''.  In formalizing this statement, we let $[uv]_T$ denote
the set of all vertices on the path in a tree $T$ that joins a vertex
$u$ and a vertex $v$.

\begin{theorem}\label{sholander}
  Let $V$ be a finite set.  A mapping $(u,v)\mapsto [uv]$ from $V^2$
  to $2^V$ satisfies {\rm (S), (T), (U$_1$)} if and only if there is a
  tree $T$ with vertex set $V$ such that $[vw]_T=[vw]$ for all pairs
  $v,w$ of its vertices. 
\end{theorem}
Sholander does not prove this theorem, but goes on to derive from the conjunction
of (S), (T), and (U$_1$) a number of corollaries that include
\begin{center}
\begin{tabular}{lp{14.5cm}}
$(2.1)$ & $\forall a,b,c\in V : b\in [a,c]\Leftrightarrow [a,b]\cap [b,c]=\{ b\}\Leftrightarrow [a,b]\cup [b,c]=[a,c]$,\\
$(5.2)$ & $\forall a,b,x,y\in V : x,y\in [a,b]\:\Rightarrow\: 
(x\in [a,y]\wedge y\in [x,b])\vee(y\in [a,x]\wedge x\in [y,b])$.
\end{tabular}
\end{center}
We are going to derive Theorem \ref{sholander} from Sholander's results.\\

The following fact is well known (for instance, Exercise 12 in Section
2.3., p.~314 of \cite{kn} and the answer on p.~558). We give its
straightforward proof just for the sake of completeness.

\begin{lemma}\label{treeorder}
Let $V$ be a finite set and let $r$ be an element of $V$. 
If $\preceq$ is a partial order on $V$ such that 
\begin{center}
\begin{tabular}{lp{14.5cm}}
{\rm (i)} & $\forall w\in V:r\preceq w$\\
{\rm (ii)} & $\forall u,v,w\in V:(u\preceq w\wedge v\preceq w)\:\Rightarrow\:(u\preceq v\vee v\preceq u)$
\end{tabular}
\end{center}
then there is a tree $T$ with vertex set $V$ such that 
$u\preceq x \;\Leftrightarrow\; u\in [rx]_T$.
\end{lemma} {\it Proof:} By induction on $|V|$. If $|V|=1$, then $T$
consists of a single vertex.  If $|V|>1$, then enumerate the minimal
elements of $V-\{r\}$ as $r_1,r_2,\ldots ,r_k$ and set $V_i=\{x\in V:
r_i\preceq x\}$. Property (ii) guarantees that the sets
$V_1,V_2,\ldots,V_k$ form a partition of $V\setminus \{ r\}$.  By the
induction hypothesis, there are trees $T_1,T_2,\ldots ,T_k$ such that
each $T_i$ has $V_i$ for its vertex set and such that elements
$u,x$ of $V_i$ satisfy $u\preceq x$ if and only if $u$ is on the path
from $r_i$ to $x$ in $T_i$.  The union of $T_1,T_2,\ldots ,T_k$ along
with vertex $r$ and the $k$ edges $rr_1,rr_2,\ldots ,rr_k$ has the
property required of $T$.  $\Box$

\medskip

\noindent {\it Proof of Theorem~\ref{sholander}.}
  The ``if'' part is clear. To prove the ``only if'' part, choose an
  arbitrary element of $V$, call it $r$, and write $u\preceq x$ if and
  only if $u\in [rx]$. This binary relation is a partial order: (1.2)
  with $a=r$ means that $\preceq$ is reflexive, (1.7) with $a=r$ means
  that $\preceq$ is antisymmetric, and (1.5) with $a=r$ implies that
  $\preceq$ is transitive. By (1.2) and (1.4) with $b=r$, $a=x$, this
  partial order has property (i) of Lemma~\ref{treeorder}; by (5.2)
  with $a=r$, $b=w$, $x=u$, $y=v$, it has property (ii) of
  Lemma~\ref{treeorder}. This lemma guarantees that there is a tree $T$
  with vertex set $V$ such that 
\begin{center}
\begin{tabular}{lp{14.5cm}}
($\alpha$) & $[rx]_T=[rx]$ for all vertices $x$ of $T$.
\end{tabular}
\end{center}
We will prove that this $T$ has the property specified in the
theorem. To begin, let us generalize ($\alpha$) to
\begin{center}
\begin{tabular}{lp{14.5cm}}
($\beta$) & if $u\in[rx]_T$, then $[ux]_T=[ux]$:
\end{tabular}
\end{center}
\noindent to verify this, note that (2.1)
with $a=r$, $b=u$, $c=x$ implies $[ru]\cap[ux]=\{u\}$ and
$[ru]\cup[ux]=[rx]$, and so 
$[ux]=([rx]\setminus [ru])\cup\{u\}
=([rx]_T\setminus [ru]_T)\cup\{u\}=[ux]_T$.
The conclusion of the theorem is a generalization of ($\beta$):
\begin{center}
\begin{tabular}{lp{14.5cm}}
($\gamma$) & $[vw]_T=[vw]$ for all pairs
  $v,w$ of vertices of $T$.
\end{tabular}
\end{center}
To verify ($\gamma$), consider an arbitrary pair $v,w$ of vertices of
$T$.  Since $T$ contains no cycle, there is a vertex $u$ such that
$[rv]_T\cap [rw]_T = [ru]_T$ and $[vw]_T=[vu]_T\cup[uw]_T$. By (1.4),
we have $[vr]\cap [rw] = [rv]\cap [rw] = [rv]_T\cap [rw]_T = [ru]_T =
[ru]$, and so (1.10) with $a=v$, $b=r$, $c=w$, $d=u$ guarantees that
$[vu]\cap[uw]=\{u\}$. Now (2.1) with $a=v$, $b=u$, $c=w$ implies
$[vu]\cup[uw]=[vw]$; using ($\beta$) and (1.4), we conclude that
$[vw]_T=[vu]_T\cup[uw]_T=[uv]_T\cup[uw]_T=[uv]\cup[uw]=[vu]\cup[uw]=[vw]$.
$\Box$

\medskip

\noindent Our proofs of Lemma~\ref{treeorder} and
Theorem~\ref{sholander} yield an efficient way of reconstructing a
tree from its collection of segments $[uv]$. Of course, the simplest
way of doing that is to make distinct $u$ and $v$ adjacent if and only
if $[uv]=\{u,v\}$.

\section{Strict Tree Betweenness and Tree Betweenness}

A {\it ternary relation\/} on a set $V$ means a subset of $V^3$; a
ternary relation ${\cal B}$ is called {\it strict} if $(x,y,z)\in
{\cal B}$ implies that $x$, $y$, and $z$ are pairwise distinct. Given
a ternary relation $\cal B$ on a set $V$, we follow Burigana~\cite{bu}
in writing $N(u,v,w)$ to mean that $u,v,w$ are pairwise distinct
elements of $V$ and $(u,v,w)\not\in {\cal B}$, $(v,w,u)\not\in {\cal
  B}$, $(w,u,v)\not\in {\cal B}$.
\begin{theorem}\label{strict}
  Let $V$ be a finite set.  A strict ternary relation $\cal B$ on $V$
  is a strict tree betweenness if and only if it satisfies 
\begin{center}
\begin{tabular}{lp{14.5cm}}
$(S_1)$ & $\forall u,v,w\in V: (u,v,w)\in {\cal B} \:\Rightarrow\: (w,v,u)\in {\cal B}$,\\
$(S_2)$ & $\forall u,v,w,z\in V: (u,v,w),(v,w,z)\in {\cal B}\:\Rightarrow\: (u,w,z)\in {\cal B}$,\\
$(S_3)$ & $\forall u,v,w,z\in V: (u,v,w),(u,w,z)\in {\cal B} \:\Rightarrow\: (v,w,z)\in {\cal B}$,\\
$(S_4)$ & $\forall u,v,w\in V: N(u,v,w) \:\Rightarrow\: \exists\, c\in V: (u,c,v),(u,c,w)\in {\cal B}$.
\end{tabular}
\end{center}
\end{theorem}
\noindent {\it Proof.\/} 
The ``only if'' part is clear. To prove the ``if'' part, we first
derive from $(S_1)$ -- $(S_4)$ a few corollaries:
\begin{center}
\begin{tabular}{lp{14.5cm}}
$(S_5)$&$\forall u,v,w,z\in V: (u,v,w),(u,w,z)\in {\cal B} \:\Rightarrow\: (u,v,z)\in {\cal B}$,\\
$(S_6)$&$\forall u,v,w,z\in V: (u,v,z),(v,w,z)\in {\cal B} \:\Rightarrow\: (u,v,w),(u,w,z)\in {\cal B}$,\\
$(S_7)$ & $\forall u,v,w\in V: (u,v,w)\in {\cal B}\:\Rightarrow\: (v,u,w)\not\in {\cal B}$,\\
$(S_8)$&$\forall u,v,w,z\in V: (u,v,z),(u,w,z)\in {\cal B} \:\Rightarrow\: 
v=w \,\vee\, (u,v,w)\in {\cal B} \,\vee\, (u,w,v)\in {\cal B}$,\\
$(S_9)$&$\forall u,v,w,z\in V: (u,v,z),(u,w,z)\in {\cal B} \:\Rightarrow\: 
v=w \,\vee\, (w,v,u)\in {\cal B} \,\vee\, (w,v,z)\in {\cal B}$,\\
$(S_{10})$&$\forall r,u,x,y,z\in V: (r,u,x),(r,u,z),(x,y,z)\in {\cal B} \:\Rightarrow\: 
y=u \,\vee\, (r,u,y)\in {\cal B}$.
\end{tabular}
\end{center}
In these derivations, we will invoke $(S_1)$ only tacitly whenever we use it.
(Whenever we invoke {\it reversed} $(S_i)$ 
we mean that we invoke the conjunction of $(S_1)$ and $(S_i)$.)

Property $(S_5)$ comes directly out of $(S_3)$ followed by $(S_2)$.
Property $(S_6)$ comes directly out of reversed $(S_3)$ followed by
$(S_2)$.  To derive $(S_7)$, note that $(w,v,w)\not\in {\cal B}$ as
${\cal B}$ is strict and that $(u,v,w)\in {\cal B}$, $(w,v,w)\not\in
{\cal B}$ implies $(w,u,v)\not\in {\cal B}$ by $(S_2)$.

We will derive $(S_8)$ and $(S_9)$ along the lines of Burigana's
proof~\cite{bu} of his Lemma 1 (i).

To derive $(S_8)$, assume the contrary: $(u,v,z),(u,w,z)\in {\cal B}$
but $(u,v,w)\not\in {\cal B}, (u,w,v)\not \in {\cal B}$ for some
$u,v,w,z$ in $V$ such that $v\ne w$. From $(u,v,z)\in {\cal B}$, we
get $(v,u,z)\not\in {\cal B}$ by $(S_7)$; in turn, from $(z,w,u)\in
{\cal B}$ and $(z,u,v)\not\in {\cal B}$, we get $(w,u,v)\not\in {\cal
  B}$ by $(S_2)$. Now $N(u,v,w)$, and so two different applications of
$(S_4)$ give points $c$ and $d$ such that $(w,c,u),(w,c,v)\in {\cal
  B}$ and $(v,d,u),(v,d,w)\in {\cal B}$. From $(u,c,w)\in {\cal B}$
and $(u,w,z)\in {\cal B}$, we get $(c,w,z)\in {\cal B}$ by $(S_3)$; in
turn, from $(v,c,w)\in {\cal B}$ and $(c,w,z)\in {\cal B}$, we get
$(v,w,z)\in {\cal B}$ by $(S_2)$. Similarly, from $(u,d,v)\in {\cal
  B}$ and $(u,v,z)\in {\cal B}$, we get $(d,v,z)\in {\cal B}$ by
$(S_3)$; in turn, from $(w,d,v)\in {\cal B}$ and $(d,v,z)\in {\cal
  B}$, we get $(w,v,z)\in {\cal B}$ by $(S_2)$. But then $(S_7)$ is
contradicted by $(v,w,z),(w,v,z)\in {\cal B}$.

Property $(S_9)$ comes out of $(S_8)$ followed by $(S_3)$ with $v$ and $w$ switched.

To derive $(S_{10})$, assume that $(r,u,x),(r,u,z),(x,y,z)\in {\cal
  B}$ and write $a\prec b$ if and only if $(r,a,b)\in{\cal B}$. This
binary relation is a strict partial order: it is irreflexive since
$\cal B$ is strict and it is transitive by $(S_5)$. By assumption, the
set $\{v\mid v\prec x, v\prec z\}$ is nonempty; consider any of its
maximal elements and denote it $w$. By $(S_8)$ and by maximality of
$w$, we have $w=u$ or $u\prec w$, and so $(S_5)$ reduces proving $y=u
\,\vee\, (r,u,y)\in {\cal B}$ to proving $y=w \,\vee\, (r,w,y)\in
{\cal B}$.
By maximality of $w$, no $c$ with $w\prec c$ satisfies $c\prec x$,
$c\prec z$; from reversed $(S_5)$, it follows that no $c$ satisfies
$(w,c,x),(w,c,z)\in {\cal B}$; since $\cal B$ is strict, $w,x,z$ are
pairwise distinct; now $(S_4)$ implies that at least one of $(w,x,z)$,
$(x,z,w)$, $(z,w,x)$ belongs to $\cal B$.  Interchangeability of $x$
and $z$ allows us to assume that at least one of $(w,x,z)$, $(z,w,x)$
belongs to $\cal B$.  In case $(w,x,z)\in {\cal B}$, we get first
$(w,x,y)\in {\cal B}$ by $(S_6)$ and then $(r,w,y)\in {\cal B}$ by
reversed $(S_2)$.  In case $(z,w,x)\in {\cal B}$, property $(S_9)$
guarantees that $y=w$ or $(w,y,x)\in {\cal B}$ or $(w,y,z)\in {\cal
  B}$; if $(w,y,x)\in {\cal B}$ or $(w,y,z)\in {\cal B}$, then
$(r,w,y)\in {\cal B}$ by reversed $(S_3)$.

Now $(S_5)$ -- $(S_{10})$ are established and we proceed to prove the
``if'' part of the theorem by induction on $|V|$. If $|V|=1$, then the
statement is trivial.  If $|V|>1$, then we choose an arbitrary element
of $V$, call it $r$, and write $a\prec b$ if and only if
$(r,a,b)\in{\cal B}$. This binary relation is a strict partial order:
it is irreflexive since $\cal B$ is strict and it is transitive by
$(S_5)$. Enumerate the minimal elements of $V-\{r\}$ as
$r_1,r_2,\ldots ,r_k$ and set $V_i=\{r_i\}\cup\{b\in V\mid r_i\prec
b\}$. Property $(S_8)$ guarantees that the sets $V_1,V_2,\ldots,V_k$
form a partition of $V\setminus \{ r\}$.  By the induction hypothesis,
there are trees $T_1,T_2,\ldots ,T_k$ such that each $T_i$ has $V_i$
for its vertex set and such that elements $x,y,z$ of each $V_i$
satisfy $(x,y,z)\in{\cal B}$ if and only if $y$ is an internal vertex
of the path in $T_i$ that joins $x$ and $z$. Let $T$ denote the union
of $T_1,T_2,\ldots ,T_k$ along with vertex $r$ and the $k$ edges
$rr_1,rr_2,\ldots ,rr_k$. We claim that elements $x,y,z$ of $V$
satisfy $(x,y,z)\in{\cal B}$ if and only if $y$ is an internal vertex
of the path in $T$ that joins $x$ and $z$. Interchangeability of $x$
and $z$ allows us to distinguish between three cases:

{\sc Case 1:} $x,z\in V_i$ {\em for some $i$.\/} In this case, the
path $P$ in $T$ that joins $x$ and $z$ is a path in $T_i$. If $y$
is an internal vertex of $P$, then $y\in V_i$ and so the induction
hypothesis guarantees that $(x,y,z)\in{\cal B}$; conversely, if
$(x,y,z)\in{\cal B}$, then $y\in V_i$ (by reversed $(S_3)$ if $r_i$
is one of $x,z$ and by $(S_{10})$ otherwise), and so the induction
hypothesis guarantees that $y$ is an internal vertex of $P$.

{\sc Case 2:} $x=r$, $z\in V_i$ {\em for some $i$.\/}

{\sc Subcase 2.1:} $z=r_i$.  In this subcase, the path in $T$ that
joins $x$ and $z$ consists of a single edge, and so it has no internal
vertex. Minimality of $r_i$ guarantees that there is no $y$ such that
$(x,y,z)\in{\cal B}$.

{\sc Subcase 2.1:} $z\ne r_i$.  In this subcase, $y$ is an internal
vertex of the path in $T$ that joins $x$ and $z$ if and only if
$y=r_i$ or $y$ is an internal vertex of the path in $T_i$ that joins
$r_i$ and $z$; our analysis of Case 1 shows that this occurs if and
only if $y=r_i$ or $(r_i,y,z)\in{\cal B}$; reversed property $(S_5)$
guarantees that $y=r_i\vee (r_i,y,z)\in{\cal B}$ implies
$(x,y,z)\in{\cal B}$; property $(S_9)$ combined with the minimality of
$r_i$ guarantees that $(x,y,z)\in{\cal B}$ implies $y=r_i\vee
(r_i,y,z)\in{\cal B}$.

{\sc Case 3:} $x\in V_i$, $z\in V_j$ {\em for some distinct $i$ and
  $j$.\/} In this case, we claim that $(x,r,z)\in{\cal B}$; to justify
this claim, let us assume the contrary.  Since $x\in V_i$ and
$(r,r_i,z)\not\in{\cal B}$, property $(S_5)$ implies
$(r,x,z)\not\in{\cal B}$; similarly, since $z\in V_j$ and
$(r,r_j,x)\not\in{\cal B}$, property $(S_5)$ implies
$(r,z,x)\not\in{\cal B}$; now $(S_4)$ gives a $c$ such that
$(r,c,x)\in{\cal B}$, $(r,c,z)\in{\cal B}$. Since $z\not\in V_i$, we
have $(r,r_i,z)\not\in{\cal B}$; in particular, $c\ne r_i$. Since
$(r,c,z)\in{\cal B}$ and $(r,r_i,z)\not\in{\cal B}$, we have
$(r,r_i,c)\not\in{\cal B}$ by $(S_5)$.  Since $(r,c,x)\in{\cal B}$,
minimality of $r_i$ implies $x\ne r_i$; in turn, $x\in V_i$ implies
$(r,r_i,x)\in{\cal B}$.  Now $(r,c,x)\in{\cal B}$, $(r,r_i,x)\in{\cal
  B}$, $c\ne r_i$, $(r,r_i,c)\not\in{\cal B}$, and so $(S_8)$ implies
$(r,c,r_i)\in{\cal B}$, contradicting minimality of $r_i$.  This
contradiction proves that $(x,r,z)\in{\cal B}$.

A vertex $y$ is an internal vertex of the path in $T$ that joins $x$
and $z$ if and only if $y=r$ or $y$ is an internal vertex of the path
in $T_i$ that joins $x$ and $r$ or $y$ is an internal vertex of the
path in $T_j$ that joins $r$ and $z$; our analysis of Case 2 shows
that this occurs if and only if $y=r$ or $(x,y,r)\in {\cal B}$ or
$(r,y,z)\in {\cal B}$; property $(S_5)$ and its reversal guarantee
that $y=r\vee (x,y,r)\in {\cal B}\vee (r,y,z)\in {\cal B}$ implies
$(x,y,z)\in{\cal B}$; property $(S_9)$ guarantees that
$(x,y,z)\in{\cal B}$ implies $y=r\vee (x,y,r)\in {\cal B}\vee
(r,y,z)\in {\cal B}$. $\Box$

\medskip

\noindent Our proof of Theorem~\ref{strict} yields an efficient way of
reconstructing a tree from its strict betweenness ${\cal B}$. Of
course, the simplest way of doing that is to make distinct $u$ and $w$
adjacent if and only if no $v$ satisfies $(u,v,w)\in {\cal B}$. 

\begin{corollary} {\rm(Burigana \cite{bu})}.
  Let $V$ be a finite set.  A strict ternary relation $\cal B$ on $V$
  is a strict tree betweenness if and only if it satisfies 
\begin{center}
\begin{tabular}{lp{14.5cm}}
$\bullet$ & $\forall u,v,w\in V: (u,v,w)\in {\cal B} \:\Rightarrow\: (w,v,u)\in {\cal B}$,\\
$\bullet$ & $\forall u,v,w\in V: (u,v,w)\in {\cal B}\:\Rightarrow\: (v,u,w)\not\in {\cal B}$,\\
$\bullet$ & $\forall u,v,w,z\in V: (u,v,w),(v,w,z)\in {\cal B}\:\Rightarrow\: (u,w,z)\in {\cal B}$,\\
$\bullet$ & $\forall u,v,w,z\in V: (u,v,w),(u,w,z)\in {\cal B} \:\Rightarrow\: (v,w,z)\in {\cal B}$,\\
$\bullet$ & $\forall u,v,w\in V: N(u,v,w) \:\Rightarrow\: \exists\, c\in V: (u,c,v),(u,c,w),(v,c,w)\in {\cal B}$.
\end{tabular}
\end{center}
$\Box$
\end{corollary}

Clearly, a ternary relation $\cal C$ on a finite set $V$ is a tree
betweenness if and only if it is the union of ternary relations $\cal
A$ and $\cal B$ such that $\cal A$ consists of all triples $(u,v,w)$
in $V^3$ that satisfy $u=v$ or $v=w$ (or both) and $\cal B$ is the
strict tree betweenness of a tree with vertex set $V$. This
observation enables us to translate our characterization of strict
tree betweenness into a characterization of tree betweenness. 

\begin{corollary}\label{weak}
  Let $V$ be a finite set.  A ternary relation $\cal C$ on $V$
  is a tree betweenness if and only if it satisfies 
\begin{center}
\begin{tabular}{lp{14.5cm}}
$(T_1)$ & $\forall u,v,w\in V: (u,v,w)\in {\cal C} \:\Rightarrow\: (w,v,u)\in {\cal C}$,\\
$(T_2)$ & $\forall u,v,w,z\in V: (u,v,w),(v,w,z)\in {\cal C}, v\ne w\:\Rightarrow\: (u,w,z)\in {\cal C}$,\\
$(T_3)$ & $\forall u,v,w,z\in V: (u,v,w),(u,w,z)\in {\cal C} \:\Rightarrow\: (v,w,z)\in {\cal C}$,\\
$(T_4)$ & $\forall u,v,w\in V: N(u,v,w) \:\Rightarrow\: \exists\, c\in V: c\ne u$ and $(u,c,v),(u,c,w)\in {\cal C}$.\\
$(T_5)$ & $\forall u,v,w\in V: (u,v,w),(v,u,w)\in {\cal C} \:\Leftrightarrow\: u=v$.\\
\end{tabular}
\end{center}
\end{corollary}
\noindent {\it Proof.\/} 
The ``only if'' part is clear. To prove the ``if'' part, assume that
$\cal C$ satisfies $(T_1)$ -- $(T_5)$, let $\cal B$ denote the set of
all triples $(u,v,w)$ in $\cal C$ such that $u,v,w$ are pairwise
distinct, and set ${\cal A}={\cal C}\setminus {\cal B}$. Clearly,
$\cal B$ satisfies $(S_1)$ -- $(S_4)$, and so it is a strict tree
betweenness. By $(T_5)$, all triples $(u,v,w)$ in $V^3$ that satisfy
$u=v$ belong to $\cal A$; in turn, by $(T_1)$, all triples $(u,v,w)$
in $V^3$ that satisfy $v=w$ belong to $\cal A$; now $(T_5)$ guarantees
that all triples $(u,v,u)$ in $\cal A$ satisfy $v=u$. 
$\Box$\\

None of the four conditions $(S_1)$ -- $(S_4)$ of Theorem~\ref{strict}
is implied by the conjunction of the other three and none of the five
conditions $(T_1)$ -- $(T_4)$ of Corollary~\ref{weak} is implied by
the conjunction of the other four. To verify this, consider
$V=\{u,v,w,z\}$ and the following five ternary relations on $V$:
\begin{eqnarray*}
{\cal B}_1 &=& \{(u,v,w),(u,v,z),(u,w,z),(v,w,z)\},\\
{\cal B}_2 &=& \{(u,v,w),(v,w,z),(w,z,u),(z,u,v),\\
&& \,\,\,
(w,v,u),(z,w,v),(u,z,w),(v,u,z)\},\\
{\cal B}_3 &=& \{(u,v,w),(u,v,z),(u,w,z),(w,v,z),\\
&& \,\,\,(w,v,u),(z,v,u),(z,w,u),(z,v,w)\},\\
{\cal B}_4 &=& \{(u,z,v),(u,z,w),(v,z,w),\\
&& \,\,\,(v,z,u),(w,z,u),(w,z,v)\},\\
{\cal B}_5 &=& V^3.
\end{eqnarray*}
For each $i=1,2,3,4$, relation ${\cal B}_i$ satisfies all seven
conditions except $(S_i)$ and $(T_i)$; relation ${\cal B}_5$ satisfies
all seven conditions except $(T_5)$.

\bigskip

\begin{center}
{\bf Acknowledgment}
\end{center}
Chv\' atal's research was undertaken, in part, thanks to funding from the Canada Research Chairs program.

\end{document}